\documentclass{bcp}

 \usepackage{amsmath,amsthm}
 \usepackage{amsfonts,amscd}
 \usepackage{amssymb,latexsym}

 \input{prepictex.tex}
 \input{pictex.tex}
 \input{postpictex.tex}

\newtheorem{theorem}{Theorem}[section]
\newtheorem{corollary}[theorem]{Corollary}

\newtheorem{proposition}[theorem]{Proposition}

\newtheorem{problem}[theorem]{Problem}
\newtheorem{remark}[theorem]{Remark}

\def\A{{\mathcal A}}

\def\D{{\mathcal D}}
\def\F{{\mathcal F}}

\def\N{{\mathcal N}}
\def\O{{\mathcal O}}
\def\P{{\mathcal P}}
\def\SS{{\mathcal S}}
\def\U{{\mathcal U}}
\def\W{{\mathcal W}}

\def\bC{{\mathbb C}}
\def\bN{{\mathbb N}}

\def\bZ{{\mathbb Z}}

\def\fG{{\mathfrak G}_n}
\def\fIG{{\mathfrak I\mathfrak G}_n}

\def\auto{{\operatorname{Aut}}}
\def\out{{\operatorname{Out}}}
\def\endo{{\operatorname{End}}}
\def\sp{{\operatorname{span}}}

\def\Ad{{\operatorname{Ad}}}
\def\id{{\operatorname{id}}}

\def\la{\langle}
\def\ra{\rangle}

\begin{document}

\mathclass{46L05, 46L40}
\abbrevauthors{R. Conti, J. H. Hong, W. Szyma{\'n}ski}
\abbrevtitle{Endomorphism of Cuntz algebras}

\title{Endomorphisms of the Cuntz Algebras}

\author{Roberto Conti~$^\dag$, 
Jeong Hee Hong~$^\ddag$
and Wojciech Szyma{\'n}ski~$^\S$
}

\address{$^\dag$~Dipartimento di Scienze, 
Universit{\`a} di Chieti-Pescara `G. D'Annunzio' \\
Viale Pindaro 42, I--65127 Pescara, Italy \\
E-mail: {\tt conti@sci.unich.it}}
\address{$^\ddag$~Department of Data Information, Korea Maritime University \\
Busan 606-791, South Korea \\
E-mail: {\tt hongjh@hhu.ac.kr}}
\address{$^\S$~Department of Mathematics and Computer Science, The University of Southern Denmark \\
Campusvej 55, DK-5230 Odense M, Denmark \\
E-mail: {\tt szymanski@imada.sdu.dk}}

\maketitlebcp

\begin{abstract}
This mainly expository article is devoted to recent advances in the study of dynamical aspects 
of the Cuntz algebras $\O_n$, $n<\infty$, via their automorphisms and, more generally, endomorphisms. 
A combinatorial description of permutative automorphisms of $\O_n$
in terms of labeled, rooted trees is presented. This in turn gives rise to 
an algebraic characterization of the restricted Weyl group of $\O_n$. 
It is shown how this group is related to certain classical dynamical systems on 
the Cantor set. An identification of the image in $\out(\O_n)$ of 
the restricted Weyl group with the group of automorphisms 
of the full two-sided $n$-shift is given, for prime $n$, providing an answer to a question raised by 
Cuntz in 1980. Furthermore, we discuss proper endomorphisms of $\O_n$ which preserve 
either the canonical UHF-subalgebra or the diagonal MASA, and present methods for constructing 
exotic examples of such endomorphisms. 
\end{abstract}

\section{Introduction}

The $C^*$-algebras $\O_n$, $n \in \{2,3,4,\ldots\} \cup \{\infty\}$ 
were first defined and investigated by Cuntz in his seminal paper 
\cite{Cun1}, and they bear his name ever since. It is difficult to overestimate the importance 
of the Cuntz algebras in theory of operator algebras and many other areas. It suffices to mention that 
Cuntz's original article, \cite{Cun1}, is probably the most cited ever paper in the area of operator 
algebras (MSC class 46L). Indeed, 
as $C^*$-algebras naturally generated by Hilbert spaces, 
the Cuntz algebras continue to provide a convenient framework for several different areas of investigations.
In order to illustrate the variety of applications, 
without pretending in any way to be exhaustive,
we only mention a very small sample of papers dealing with Fredholm theory, 
classification of $C^*$-algebras,
self-similar sets, coding theory, continuous fractions, spectral flow and index theory for 
twisted cyclic cocycles, see e.g. \cite{Bi,KP,MSW,MM,KHL,CPR}.

This mainly expository article is devoted to recent advances in the study of dynamical aspects 
of the Cuntz algebras $\O_n$ with $n<\infty$ via their automorphisms and, 
more generally, endomorphisms. It is not 
a comprehensive review but a selective one, biased towards the contributions made by the 
three authors. Some original results are also presented in this article, as will be explained later. 

Systematic investigations of endomorphisms of $\O_n$, $n<\infty$ 
were initiated by Cuntz in \cite{Cun2}.   
A fundamental bijective correspondence between unital $*$-endomorphisms and unitaries in $\O_n$ 
was established therein (see equation (\ref{Cuntzendo}), below). Using this correspondence Cuntz 
proved a number of interesting results, in particular with regard to those endomorphisms which 
globally preserve either the core UHF-subalgebra $\F_n$ or the diagonal MASA $\D_n$. 

Investigations of automorphisms of $\O_n$ began almost immediately after the birth 
of the algebras in question, see \cite{A,Cun2,ETW,EFTW,CaEv,MT,T}. Classification of group actions 
on $\O_n$ came to the fore somewhat later, see \cite{I3,Matui}. In the present article, 
we review more recent results on automorphisms of $\O_n$ contained in \cite{Szy,CS,CKS1,CKS2,CHS}. 

Proper endomorphisms of the Cuntz algebras have also attracted a lot of attention. 
In particular, they played a role in certain aspects of  index 
theory, both from the $C^*$-algebraic and von Neumann algebraic point of view. 
The problem of computing the Jones(-Kosaki-Longo) index of (the normal extensions) 
of localized endomorphisms of $\O_n$ was posed in \cite{J}. Progress on this and other 
related problems was then achieved in a number of papers. Of particular note 
in this regard are contributions made by Longo, \cite{Lo1,Lo2,Lo3}, and Izumi, \cite{I0,I1,I2}, 
but see also \cite{Cun3,CP,Ak,CF,G,Kaw1,Kaw2,CS-ind,CRS,HSS}. 
There is also a parallel line of reasearch dealing with various entropy 
computations, e.g. see \cite{Choda,Sk,SZ}.
Recently, one of the most interesting applications of 
endomorphisms of $\O_n$, found by Bratteli and J{\o}rgensen in \cite{BJ1,BJ2}, is in the 
area of wavelets. Before that, shift endomorphisms of Cuntz algebras have been systematically
employed in the analysis of structural aspects of quantum field theory, 
see e.g. the discussion in \cite[Section 2]{DR} and references therein.

The present article is organized as follows. After setting the stage with some preliminaries in Section 2, 
we discuss localized automorphisms in Section 3. Localization refers to the fact that the corresponding 
unitary lies in one of the matrix algebras constituting a building block of the UHF-subalgebra $\F_n$. 
In section \ref{palt}, we review fundamental results about permutative automorphisms of $\O_n$, mainly contained 
in \cite{CS}. The key breakthrough obtained therein was a clear-cut correspondence between 
such automorphisms and certain combinatorial structure related to labeled trees. This in turn  
served as a platform for further theoretical analysis, classification results, and construction of non-trivial examples.

In Section \ref{pme}, we present a
more direct approach to finding automorphisms, based on solving certain polynomial matrix
equations. Even though these equations are relatively easy to derive, finding a complete
set of solutions is a highly non-trivial task.

Section \ref{class} contains a complete classification of those permutative endomorphisms of $\O_3$ in 
level $k=3$ which are either automorphisms of $\O_3$ or restrict to automorphisms of the diagonal 
$\D_3$. These results were obtained in \cite{CKS1} and in the subsequent unpublished work \cite{CKS2}, 
with aid of massive computer calculations. We also give tables summarizing
the results of our automorphism search for all values of parameters $n$ and $k$ with $n+k \leq 6$.

In Section \ref{rweyl}, we review very recently obtained description of the so-called restricted  
Weyl group of $\O_n$ in terms of automorphisms of the full two-sided $n$-shift,  \cite{CHS}. On one hand, 
this result provides an answer to a question raised by Cuntz in \cite{Cun2}. On the other hand, 
it establishes a very interesting correspondence between an important class of automorphisms of a 
purely infinite, simple $C^*$-algebra $\O_n$ and much studied automorphism group of a classical system 
of paramount importance in symbolic dynamics, \cite{Kit,LM}. Some aspects of this correspondence 
are related to the problem of extension of an automorphism from $\D_n$ to the entire $\O_n$. 
A similar question for the UHF-subalgebra $\F_n$ rather than the diagonal $\D_n$ was studied 
recently in \cite{Conti}. 

Some recent results related to proper endomorphisms of $\O_n$ are reviewed in Section \ref{propend}. 
Subsection \ref{endpresUHF} deals with those endomorphisms which globally preserve the UHF-subalgebra $\F_n$, 
while Subsection \ref{endpresdiag} with those which globally preserve the diagonal $\D_n$. 
The main theme in here is construction of endomorphisms which globally preserve one of these subalgebras 
but whose corresponding unitary does not belong to the relevant normalizer. The problem of existence of such 
exotic endomorphisms was left open in \cite{Cun2} and remained unresolved until the recent works 
of \cite{CRS} and \cite{HSS}. 

\section{Preliminaries}

If $n$ is an integer greater than 1, the Cuntz algebra $\O_n$ is the 
unital $C^*$-algebra generated by $n$ isometries $S_1, \ldots, S_n$, satisfying
$\sum_{i=1}^n S_i S_i^* = I$, \cite{Cun1}.
Then it turns out that $\O_n$ is
separable, simple, nuclear and purely infinite.
We denote by $W_n^k$ the set of $k$-tuples $\mu = (\mu_1,\ldots,\mu_k)$
with $\mu_m \in \{1,\ldots,n\}$, and by $W_n$ the union $\cup_{k=0}^\infty W_n^k$,
where $W_n^0 = \{0\}$. We call elements of $W_n$ multi-indices.
If $\mu \in W_n^k$ then $|\mu| = k$ is the length of $\mu$.
If $\mu = (\mu_1,\ldots,\mu_k) \in W_n$, then $S_\mu = S_{\mu_1} \ldots S_{\mu_k}$
($S_0 = 1$ by convention) is an isometry with range projection $P_\mu=S_\mu S_\mu^*$.
Every word in $\{S_i, S_i^* \ | \ i = 1,\ldots,n\}$ can be uniquely expressed as
$S_\mu S_\nu^*$, for $\mu, \nu \in W_n$ \cite[Lemma 1.3]{Cun1}.

We denote by $\F_n^k$ the $C^*$-subalgebra of $\O_n$ spanned by all words of the form
$S_\mu S_\nu^*$, $\mu, \nu \in W_n^k$, which is isomorphic to the
matrix algebra $M_{n^k}({\mathbb C})$. The norm closure $\F_n$ of
$\cup_{k=0}^\infty \F_n^k$, is the UHF-algebra of type $n^\infty$,
called the core UHF-subalgebra of $\O_n$, \cite{Cun1}. It is the fixed point algebra
for the gauge action of the circle group $\gamma:U(1)\rightarrow{\auto}(\O_n)$ defined
on generators as $\gamma_t(S_i)=tS_i$. For $k\in\bZ$, we denote by $\O_n^{(k)}
:=\{x\in\O_n:\gamma_t(x)=t^k x\}$,
the spectral subspace for this action. In particular, $\F_n=\O_n^{(0)}$.
The $C^*$-subalgebra of $\F_n$ generated by projections $P_\mu$, $\mu\in W_n$, is a 
MASA (maximal abelian subalgebra) both in $\F_n$ and in $\O_n$. We call it
the diagonal and denote $\D_n$. The spectrum of $\D_n$ is naturally identified with 
$X_n$ --- the full one-sided $n$-shift space. We also set $\D_n^k:=\D_n\cap\F_n^k$. 
Throughout this paper we are interested in the inclusions
$$ \D_n \subseteq \F_n\subseteq \O_n.$$
The UHF-subalgebra $\F_n$ posseses a unique normalized trace, denoted $\tau$. 
We will refer to the restriction of $\tau$ to $\D_n$ as to the canonical trace on $\D_n$. 

We denote by $\SS_n$ the group of those unitaries in $\O_n$ which can be written
as finite sums of words, i.e., in the form $u = \sum_{j=1}^m S_{\mu_j}S_{\nu_j}^*$
for some $\mu_j, \nu_j \in W_n$. 
It turns out that $\SS_n$ is isomorphic to the Higman-Thompson group $G_{n,1}$ \cite{N}.
One can also identify a copy of Thompson's group $F$ sitting in canonical fashion inside $\SS_2$.
We also denote $\P_n=\SS_n\cap\U(\F_n)$. Then
$\P_n=\cup_k\P_n^k$, where $\P_n^k$ are permutation unitaries in $\U(\F_n^k)$.
That is, for each $u\in\P_n^k$ there is a unique permutation $\sigma$ of multi-indices
$W_n^k$ such that
\begin{equation}\label{permutunitary}
u = \sum_{\mu \in W_n^k} S_{\sigma(\mu)} S_\mu^*.
\end{equation}

As shown by Cuntz in \cite{Cun2}, there exists the following bijective correspondence
between unitaries in $\O_n$ and unital $*$-endomorphisms of $\O_n$ (whose collection we denote
by $\endo(\O_n)$). A unitary $u$  in $\O_n$ determines an endomorphism $\lambda_u$ 
by\footnote{In some papers, e.g. \cite{Cun2}, \cite{Szy} and \cite{CS}, a different 
convention $\lambda_u(S_i)=u^*S_i$ is used.}
\begin{equation}\label{Cuntzendo}
\lambda_u(S_i) = u S_i, \;\;\; i=1,\ldots, n. 
\end{equation}
Conversely, if $\rho :\O_n\rightarrow \O_n$ is an endomorphism, then
$\sum_{i=1}^n\rho(S_i)S_i^*=u$ gives a unitary $u\in\O_n$
such that $\rho=\lambda_u$. If the unitary $u$ arises from a permutation $\sigma$ via the formula
\eqref{permutunitary}, the corresponding endomorphism will be 
sometimes denoted by $\lambda_{\sigma}$.
Composition of endomorphisms corresponds to a `convolution'
multiplication of unitaries as follows:
\begin{equation}\label{convolution}
\lambda_u \circ \lambda_w = \lambda_{u*w}, \;\;\; \text{where} \;\; u*w=\lambda_u(w)u. 
\end{equation}
We denote by $\varphi$ the canonical shift:
$$ \varphi(x)=\sum_{i=1}^nS_ixS_i^*, \;\;\; x\in\O_n. $$
If we take $u=\sum_{i, j}S_iS_jS_i^*S_j^*$ then $\varphi=\lambda_u$.
It is well-known that $\varphi$ leaves invariant both $\F_n$ and $\D_n$, and that $\varphi$
commutes with the gauge action $\gamma$. We denote by $\phi$ the standard left inverse 
of $\varphi$, defined as $\phi(a)=\frac{1}{n}\sum_{i=1}^n S_i^* a S_i$. 

If $u\in\U(\O_n)$ then for each positive integer $k$ we denote
\begin{equation}\label{uk}
u_k := u \varphi(u) \cdots \varphi^{k-1}(u).
\end{equation}
We agree that $u_k^*$ stands for $(u_k)^*$. If
$\alpha$ and $\beta$ are multi-indices of length $k$ and $m$, respectively, then
$\lambda_u(S_\alpha S_\beta^*)=u_kS_\alpha S_\beta^*u_m^*$. This is established through
a repeated application of the identity $S_i a = \varphi(a)S_i$, valid for all
$i=1,\ldots,n$ and $a \in \O_n$. If $u\in\F_n^k$ for some $k$ then, 
following \cite{CP},  we call endomorphism $\lambda_u$ {\em localized}. Even though systematic 
investigations of such endomorphisms were initiated in \cite{Cun2}, it should be noted 
that automorphisms constructed this way appeared already in the work of Connes in the 
context of the hyperfinite type $II_1$ factor, \cite{Connes}. 

For algebras $A\subseteq B$ we denote by $\N_B(A)=\{u\in\U(B):uAu^*=A\}$ the normalizer
of $A$ in $B$, and by $A' \cap B=\{b \in B: (\forall a \in A) \; ab=ba\}$ the 
relative commutant of $A$ in $B$. We also denote by $\auto(B,A)$ the collection of all those 
automorphisms $\alpha$ of $B$ such that $\alpha(A)=A$, and by $\auto_A(B)$ 
those 
automorphisms of $B$ which fix $A$ point-wise. 

\section{Localized endomorphisms and automorphisms}

In this section, we mostly deal with automorphisms of $\O_n$. However, it may be
useful to broaden our horizon for a little while and consider more general 
endomorphisms of Cuntz algebras from the point of view of subfactor/sector theory.

\subsection{One example}\label{expls}
Since dealing simultaneously with all unitaries in matrix algebras is very difficult,
in order to discuss interesting cases
it is convenient to focus on some selected classes of unitaries which arise in specific situations
like in the study of integrable systems.
Let $H$ be a Hilbert space with $\dim(H)=n$.
Cuntz already noticed that unitary solutions $Y \in \U(H \otimes H)$ of the quantum YBE 
(without spectral parameter)
\begin{equation}\label{YBE}
Y_{12}Y_{23}Y_{12}=Y_{23}Y_{12}Y_{23}
\end{equation}
can be characterized in Cuntz algebra terms as those unitaries $Y$ in $\F_n^2$ satisfying
\begin{equation}\label{CuntzYBE}
\lambda_Y(Y) = \varphi(Y) \ . 
\end{equation}
As a simple exercise,
it is instructive to observe that no nontrivial unitary solution of the YBE
induces an automorphism of $\O_n$. Indeed, we claim
that if $Y \in \F_n^2$ then  $Y$ satisfies equation (\ref{YBE})
if and only if\footnote{At first sight 
this condition might look a bit strange, however one should then remember that in 
algebraic quantum field theory the canonical braiding $\epsilon_\rho$ of a 
localized morphism $\rho$ of the observable net $\mathfrak A$ indeed satisfies
$\epsilon_\rho \in \rho^2({\mathfrak A})' \cap {\mathfrak A}$.} 
$Y \in (\lambda_Y^2,\lambda_Y^2):=\lambda^2_Y(\O_n)' \cap \O_n$.
Here one needs the composition rule of endomorphisms, 
namely $\lambda_Y^2 = \lambda_{\lambda_Y(Y)Y}=\lambda_{Y\varphi(Y)Y\varphi(Y^*)}$,
along with the characterization of self-intertwiners recalled in Section \ref{endpresUHF} below. 
That is, 
thanks to equation (\ref{relcomm}) one has that
$Y \in (\lambda_Y^2,\lambda_Y^2)$ if and only if 
\begin{equation}\label{ybeq}
\big(Y\varphi(Y)Y\varphi(Y^*)\big)^* Y Y\varphi(Y)Y\varphi(Y^*) = \varphi(Y) \ .
\end{equation}
Now, the l.h.s. of (\ref{ybeq}) is precisely $\varphi(Y)Y^*\varphi(Y^*)Y \varphi(Y) Y \varphi(Y^*)$ 
and the claim is now clear. If $Y$ is not a multiple of the identity, 
this shows already that $\lambda_Y^2(\O_n)' \cap \O_n$
contains non-scalar elements and therefore, $\O_n$ being simple, 
$\lambda_Y^2$ is not an automorphism, as well as $\lambda_Y$.
The computation of the Jones index for subfactors associated to 
Yang-Baxter unitaries has been discussed in more detail in \cite{CP,CF}.

It is well-known that finding all solutions of the YBE in dimension $n$ is a difficult problem
that has been dealt with only for very small values of $n$.
This is closely related with the classification problem for braiding in categories of representations of
quantum groups and/or conformal nets.
It is expected that attaching to these solutions invariants from subfactor theory will lead to 
a much better understanding.

Other families of unitaries related to the study of spin/vertex models might also provide a useful playground:

\begin{problem}

\begin{itemize}
\item[(a)] Examine Cuntz algebra endomorphisms associated to normalized Hadamard matrices, cf. \cite{J};
\item[(b)] Discuss from the Cuntz algebra point of view the tetrahedron equation and/or its several variations 
(see e.g. \cite{BMS}).
\end{itemize}
\end{problem}

Finally, it is worth to recall that a throughout discussion of localized endomorphisms associated to (finite-dimensional) 
unitaries satisfying the so-called pentagon equation 
(which is a basic ingredient of quantum group theory) 
has been provided in \cite{Lo3,CP}.

\subsection{Permutative automorphisms and labeled trees}\label{palt}

We begin by recalling\footnote{Note the difference in convention regarding the definition of $\lambda_u$.} 
some results from \cite{CS}. Let $u$ be a unitary in  $\F_n^k$. 
For $i,j \in \{1,\ldots,n\}$, one defines linear maps $a_{ij}^u: \F_n^{k-1} \to \F_n^{k-1}$ by
$a_{ij}^u(x) = S_i^* u^* x u S_j$, $x \in \F_n^{k-1}$. We denote $V_u = \F_n^{k-1}/\bC 1$.
Since $a_{ij}^u(\bC1 ) \subseteq \bC 1$, there are induced maps $\tilde{a}_{ij}^u: V_u \to V_u$.
We define $A_u$ as the subring of ${\mathcal L}(V_u)$ generated by
$\{\tilde{a}_{ij}^u \ | \ i,j = 1,\ldots, n\}$. We denote by $H$ the linear span of the $S_i$'s.
Following \cite{CP}, we define inductively $\Xi_0 = \F_n^{k-1}$ and 
$\Xi_{r} = \lambda_u(H)^* \Xi_{r-1} \lambda_u(H)$, $r \geq 1$. 
It follows that $\{\Xi_r\}$ is a nonincreasing sequence of subspaces of $\F_n^{k-1}$ 
and thus it eventually stabilizes. If $p$ is the smallest integer for which $\Xi_p = \Xi_{p+1}$, then 
$\Xi_u := \bigcap_r \Xi_r = \Xi_p$. The following result is contained in \cite{CS}. 

\begin{theorem}\label{main}
Let $u$ be a unitary in $\F_n^k$. Then the following conditions are equivalent:
\begin{itemize}
\item[(1)] $\lambda_u$ is invertible with localized inverse;
\item[(2)] $A_u$ is nilpotent;
\item[(3)] $\Xi_u = {\mathbb C}1$.
\end{itemize}
\end{theorem}

In the case of a permutation unitary $u\in\P_n^k$, Theorem \ref{main} may be strengthened and 
very conveniently reformulated in combinatorial terms, as follows. As shown in \cite{CS}, the corresponding 
$\lambda_u$ is an automorphism of $\O_n$ if and only if $u$ satisfies two conditions, called $(b)$ and 
$(d)$ therein. Condition $(b)$ by itself guarantees that endomorphism $\lambda_u$ restricts 
to an automorphism of the diagonal $\D_n$.\footnote{Since $\P_n^k \subset \N_{\O_n}(\D_n)$ for all $k$,
every permutative endomorphism of $\O_n$ maps $\D_n$ into itself.} To describe these two conditions we will  identify 
unitary $u\in\P_n^k$ with the corresponding permutation of $W_n^k$. 

For $i=1,\ldots,n$, one defines a mapping $f_i^u:W_n^{k-1}\to W_n^{k-1}$ so that 
$f_i^u(\alpha)=\beta$ if and only if there exists $m\in W_n^1$ such that $(\beta,m)=u(i,\alpha)$. 
Then $u$ satisfies condition $(b)$ if and only if there exists a partial order $\leq$ on 
$W_n^{k-1} \times W_n^{k-1}$ such that:
\begin{itemize}
\item[(i)] Each element of the diagonal $(\alpha,\alpha)$ is minimal;
\item[(ii)] Each $(\alpha,\beta)$ is bounded below by some diagonal element;
\item[(iii)] For every $i$ and all $(\alpha,\beta)$ such that $\alpha \neq \beta$, we have
\begin{equation}
(f_i^u(\alpha),f_i^u(\beta)) \leq (\alpha,\beta) \ .
\end{equation}
\end{itemize}
For this condition (b) to hold it is necessary that the diagram of each mapping $f_i^u$ is a rooted 
tree\footnote{
Trees continue to be used in a number of different contexts, sometimes related 
to operator algebras, see e.g. \cite{Ar,B,LoRo,CoKr,HoOl}, 
however our approach seems to be genuinely new.}, 
with the root its unique fixed point and with an edge going down from $\alpha$ to $\beta$ if 
$f_i^u(\alpha)=\beta$. By convention, we do not include in the diagram
the loop from the root to itself. For example, if $u=\id$ is viewed as an element of $\P_2^3$,  
then the corresponding pair of labeled trees is: 
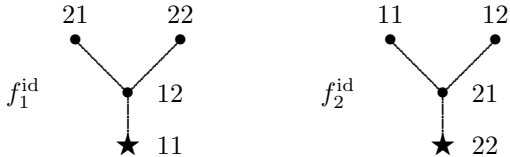
\[ \beginpicture
\setcoordinatesystem units <0.7cm,0.7cm>
\setplotarea x from 4 to 9, y from -1 to 1

\put {$\bullet$} at -1 1
\put {$\bullet$} at 1 1
\put {$\bullet$} at 0 0
\put {$\bigstar$} at 0 -1

\setlinear
\plot -1 1 0 0 /
\plot 1 1 0 0 /
\plot 0 0 0 -1 /
\put {$f_1^{\id}$} at -2 0
\put {$21$} at -1 1.5
\put {$22$} at 1 1.5
\put {$12$} at 0.8 0
\put {$11$} at 0.8 -1

\put {$\bullet$} at 5 1
\put {$\bullet$} at 7 1
\put {$\bullet$} at 6 0
\put {$\bigstar$} at 6 -1

\setlinear
\plot 5 1 6 0 /
\plot 7 1 6 0 /
\plot 6 0 6 -1 /
\put {$f_2^{\id}$} at 4 0
\put {$11$} at 5 1.5
\put {$12$} at 7 1.5
\put {$21$} at 6.8 0
\put {$22$} at 6.8 -1
\endpicture \]

To describe condition $(d)$ we define  $\W_n^{k-1}$ as the union of all off-diagonal elements 
of $W_n^{k-1}\times W_n^{k-1}$ and one additional element $\dag$. For $i,j\in W_n^1$, we also define 
mappings $f_{ij}^u:\W_n^{k-1}\to\W_n^{k-1}$ so that $f_{ij}^u(\alpha,\beta)=(\gamma,\delta)$ 
if there exists an $m\in W_n^1$ such that $(\gamma,m)=u(i,\alpha)$ and $(\delta,m)=u(j,\beta)$. 
Otherwise, we set $f_{ij}^u(\alpha,\beta)=\dag$. We also put $f_{ij}^u(\dag)=\dag$ for 
all $i,j$. Then $u$ satisfies condition $(d)$ if and only if 
there exists a partial order $\leq$ on $\W_n^{k-1}$ such that:
\begin{itemize}
\item[(i)]
The only minimal element with respect to $\leq$ is $\dagger$.
\item[(ii)]
For every $(\alpha,\beta) \in \W_n^{k-1}$ and all $i,j = 1, \ldots, n$, we have 
\begin{equation}\label{fij}
f_{ij}^u(\alpha,\beta) \leq (\alpha,\beta).
\end{equation}
\end{itemize}

With help of this combinatorial approach, a complete classification has been achieved in \cite{CS}, 
\cite{CKS1} and \cite{CKS2} of permutations in $\P_n^k$ with $n+k\leq 6$ such that the 
corresponding endomorphism $\lambda_u$ is either automorphism of $\O_n$ or restricts to 
an automorphism of the diagonal $\D_n$. Considering the image of $\lambda(\P_n^k)^{-1}$ in the 
outer automorphism group of $\O_n$, it was shown in \cite{CS} with respect to the case of $\O_2$, 
that no outer automorphisms apart from the flip-flop arise in this way for 
$k=3$ (a much simpler case $k=2$ being already known). For $k=4$, twelve new classes in 
$\out(\O_2)$ were found. 

\subsection{Inverse pairs of localized automorphisms}\label{pme}

In this section, we gather together a few facts about pairs of
unitaries in some finite matrix algebras giving rise to
automorphisms of $\O_n$ that are inverses of each other.
We also briefly discuss interesting algebraic equations such
unitaries must satisfy. These equations provide a 
useful background for the considerations in Section 3 of \cite{CS}
(e.g. Theorem 3.2, Corollary 3.3 therein), which are reviewed in the present article in 
Section \ref{palt} above.
They have also been useful for several other concrete computations in \cite{CS}, 
e.g. in computing explicitly the inverse of $\lambda_A$, introduced and analyzed 
in Section 5, filling the tables of Section 6,
and in the search of square-free automorphisms, \cite{CS}.
Although these equations are not difficult to derive, we think that highlighting them may be 
of benefit, especially to the readers who do not use the machinery of Cuntz algebras 
on the daily basis. 

%

So let us suppose that $u \in \F_n^k$ and $w \in \F_n^h$ are unitaries such that
$$\lambda_u \lambda_w = {\rm id} = \lambda_w \lambda_u  \ , $$
i.e. $\lambda_u(w)u = 1 = \lambda_w(u)w$.\footnote{Since $\lambda_u$ 
and $\lambda_w$ are injective, one identity implies the other.
Also, up to replacing $k$ and $h$ with
their maximum,
there would be no loss of generality in assuming that $k = h$. However as
the inverse of an automorphism induced by a unitary in a
matrix algebra might very well be induced by a unitary in a
larger matrix algebra, it seems convenient to allow this more flexible asymmetric formulation.
It is worth stressing that, given $k$, the subset of unitaries $u$'s in $\F_n^k$
such that $\lambda_u^{-1}$ (exists and) is still induced by a unitary in $\F_n^k$
is strictly smaller than the set of unitaries such that $\lambda_u^{-1}$ is induced
by a unitary in some $\F_n^h$.
An a priori bound for $h$ as a function of $n,k$ is provided in \cite[Corollary 3.3]{CS}.}
This readily leads to a system of coupled matrix equations 
\begin{equation}\label{coupled}
u_h w u_h^* = u^*, \quad w_k u w_k^* = w^*,
\end{equation}
where both $u_h$ and $w_k$ are in $\F_n^{h + k-1}$.
In passing, observe that the second equation is independent of the level $h$
for which $w\in\F_n^h$.

In practical situations, one is faced with the converse problem.
Starting with some $u \in \F_n^k$, one might not know the precise value of $h$,
or even if the corresponding $w$ exists at all.
It turns out that 
the existence of
solutions (for $w$) of equations (\ref{coupled}) imply
invertibility of $\lambda_u$. The following proposition combined with \cite[Corollary 3.3]{CS}
gives an algorithmic procedure for finding these solutions. We omit its elementary proof. 
\begin{proposition}
Let $u$ be a unitary in $\F_n^k$ and suppose that $u_h^* u^* u_h \in \F_n^h$ for some $h$.
Then $\lambda_u$ is invertible and $\lambda_u^{-1}=\lambda_w$ with $w := u_h^* u^* u_h$.
\end{proposition}

In particular, given a unitary $u \in \F_n^k$, one has $\lambda_u^2 = {\rm id}$
(i.e., $u = w$)
if and only if $\lambda_u(u)u = 1$, if and only if $u_k u u_k^* = u^*$.

\medskip
Finally, we present yet another computational strategy for determining invertibility of
endomorphism $\lambda_u$ and finding its inverse. Again, we omit an elementary proof
of the following proposition.
\begin{proposition}\label{equ}
Let $u$ and $w$ be unitaries in $\F_n^k$ and $\F_n^h$, respectively, satisfying
equations (\ref{coupled}). Then $u$ is a solution of the following polynomial matrix equation
\begin{equation}\label{necU}
(u_r^* u^* u_r)_r u (u_r^* u^* u_r)_r^* = u_r^* u u_r \ ,
\end{equation}
where $r$ can be taken as maximum of $k$ and $h$.

Conversely, given $r$, every solution $u \in \F_n^r$ of equation (\ref{necU})
gives rise to an automorphism $\lambda_u$ of $\O_n$,
with inverse induced by $w:=u_r^* u^* u_r$.
\end{proposition}
After some simplification, taking into account that $u \in \F_n^r$,
it is straighforward to check that the first nontrivial equation
in the family (\ref{necU}), for $r=2$, is
\begin{equation}
\varphi(u)\varphi^2(u^*)\varphi(u^*)u = u \varphi(u)\varphi^2(u^*)\varphi(u^*) \ ,
\end{equation}
i.e. $u$ commutes with $\varphi(u\varphi(u^*)u^*)$.
Similarly, for $r=3$, one obtains
\begin{eqnarray*}
&u\varphi^2(u\varphi(u)\varphi^2(u))
\varphi\Big(\varphi^2(u^*)\varphi(u^*)u\varphi(u)\varphi^2(u)\Big)\varphi^2(u^*)\varphi(u^*)\\
& =\varphi^2(u\varphi(u)\varphi^2(u))
\varphi\Big(\varphi^2(u^*)\varphi(u^*)u\varphi(u)\varphi^2(u)\Big)\varphi^2(u^*)\varphi(u^*)u
\end{eqnarray*}
i.e., $u$ commutes with $\varphi^2(u\varphi(u)\varphi^2(u))
\varphi\Big(\varphi^2(u^*)\varphi(u^*)u\varphi(u)\varphi^2(u)\Big)\varphi^2(u^*)\varphi(u^*)$.

\begin{remark}
\rm The strategy of 
applying Proposition \ref{equ} is to find all pairs satisfying (\ref{coupled})
by solving equations of the form (\ref{necU}) for all values of $r$.
Implicitly, by solving such an equation, we predict $w$ to take a
particular form, namely $w=u_r^* u^* u_r$. However, we do not assume $w\in\F_n^r$.
In fact, $w$ automatically belongs to $\F_n^{2r-1}$. Combining this with equations (\ref{coupled})
we obtain an additional relation $u$ must satisfy, namely
$u_r^* u^* u_r = u_{2r-1}^* u^* u_{2r-1}$.
\end{remark}

We find it rather intriguing that in the case of permutation unitaries
the polynomial matrix equations (\ref{necU}) turn out
to be equivalent to the tree related conditions of \cite[Corollary 4.12]{CS}.

\medskip
Of course, the above polynomial matrix equations apply to
arbitrary unitaries in the algebraic part of ${\mathcal F}_n$ and not only to
permutation matrices. Therefore, they can be used for finding other families of
automorphisms of $\O_n$ with localized inverses. It is to be expected that new
interesting classes of automorphisms different from the much studied quasi-free ones
will be found this way. It seems also worth while to investigate the algebraic varieties 
in ${\mathbb R}^{2k^2}$ 
defined by these equations. At present, we are not aware of occurences
of these equations outside the realm of Cuntz algebras but we would not be surprised
if such instances were found.

\subsection{The classification of permutative automorphisms}\label{class}
The classification of automorphisms of $\O_n$ associated to unitaries in $\P_n^k$ for $n=k=2$
goes back to \cite{Kaw1}.
Beyond that,
the program of classifying permutative automorphisms corresponding to unitaries in $\P_n^k$ for 
small values of $n$ and $k$ was 
initiated in \cite{CS} and continued in \cite{CKS1}. In this section, 
we present a complete classification in the case $n=k=3$. These results come from the 
unpublished manuscript, \cite{CKS2}, and were obtained with aid of a massive scale computer 
calculations involving Magma software, \cite{Magma}. 

As discussed in Section \ref{palt} above, determination of invertibility 
of a permutative 
endomorphism hinges upon verification of two combinatorial conditions, called $(b)$ and $(d)$, \cite{CS}. 
In short, condition $(b)$ allows to determine when the corresponding endomorphism 
$\lambda_\sigma$ of $\O_n$ restricts to an automorphism of $\D_n$, while condition $(d)$, 
together with $(b)$, determines the more stringent
situation that $\lambda_\sigma \in \auto(\O_n)$.
\footnote{It is also useful to observe that, 
the diagonal $\D_n$ being a MASA in $\O_n$,
an automorphism of $\O_n$ mapping $\D_n$ into itself automatically restricts to an automorphism of $\D_n$.}
Detailed analysis of conditions $(b)$ and $(d)$ 
in terms of labeled, rooted trees, 
was then accomplished for $n=2$ in \cite{CS} up to level $k=4$,
and in \cite{CKS1} for $n=3$ up to level $k=3$ (for $n=3=k$ only condition $(b)$ 
was examined) and $n=4$ up to level 2.

In the case $n=k=3$,
the involved rooted trees have nine vertices. 
By in-degree type of a rooted tree we mean the multiset of the in-degrees of its vertices;
in \cite[Figure 1]{CKS1}, we have divided a relevant subset of 171 rooted trees with 9 vertices into 11 
distinct in-degree types called $A,B,\ldots,K$ and described in Table 1 therein.
For instance, the in-degree type $A$ spots only trees with six vertices with no incoming edge (leaves)
and three vertices with three incoming edges 
(also recall that there is always an invisible loop at the root),
while the in-degree type $B$ singles out trees with five leaves, one vertex with one incoming edge,
one vertex with two incoming edges and two vertices with three incoming edges.
It turns out that condition $(b)$ is satisfied for a set $\F$ 
of 7390 3-tuples of labelled rooted trees,
up to permutation of tree position (action of the symmetric group $S_3$) and 
consistent relabelling of all trees (action of $S_9$),
as described in \cite[Section 2.2]{CKS1}.
The set $\F$ is then partitioned into 6 distinct three-element multisets of in-degree types,
as listed in Table \ref{some_sat_condD} below
(based on Table 2 in \cite{CKS1}, to which we refer for more details).
Examples of triples of rooted trees with labels belonging to the in-degree types
$A \, A \, A$ and $A \, F\, G$ are shown in Figure 2 of \cite{CKS1}.

For instance, the six permutative Bogolubov automorphisms associated to permutations $u \in \P_3^1$, 
viewed as elements in $\P_3^3$, give raise to the 3-tuple of trees with in-degree type $A \, A \, A$
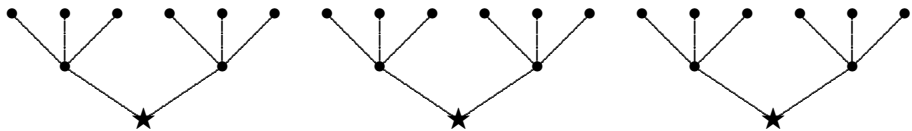
\[ \beginpicture
\setcoordinatesystem units <0.7cm,0.7cm>
\setplotarea x from 10 to 15, y from -1 to 1

\put {$\bullet$} at 0.5 1
\put {$\bullet$} at 1.5 1
\put {$\bullet$} at 2.5 1
\put {$\bullet$} at -0.5 1
\put {$\bullet$} at -1.5 1
\put {$\bullet$} at -2.5 1
\put {$\bullet$} at -1.5 0
\put {$\bullet$} at 1.5 0
\put {$\bigstar$} at 0 -1

\setlinear
\plot 0 -1 -1.5 0 /
\plot 0 -1 1.5 0 /
\plot -1.5 0 -2.5 1 /
\plot -1.5 0 -1.5 1 /
\plot -1.5 0 -0.5 1 /
\plot 1.5 0 2.5 1 /
\plot 1.5 0 1.5 1 /
\plot 1.5 0 0.5 1 /

\put {$\bullet$} at 6.5 1
\put {$\bullet$} at 7.5 1
\put {$\bullet$} at 8.5 1
\put {$\bullet$} at 3.5 1
\put {$\bullet$} at 4.5 1
\put {$\bullet$} at 5.5 1
\put {$\bullet$} at 4.5 0
\put {$\bullet$} at 7.5 0
\put {$\bigstar$} at 6 -1

\setlinear
\plot 6 -1 4.5 0 /
\plot 6 -1 7.5 0 /
\plot 4.5 0 3.5 1 /
\plot 4.5 0 4.5 1 /
\plot 4.5 0 5.5 1 /
\plot 7.5 0 6.5 1 /
\plot 7.5 0 7.5 1 /
\plot 7.5 0 8.5 1 /

\put {$\bullet$} at 12.5 1
\put {$\bullet$} at 13.5 1
\put {$\bullet$} at 14.5 1
\put {$\bullet$} at 9.5 1
\put {$\bullet$} at 10.5 1
\put {$\bullet$} at 11.5 1
\put {$\bullet$} at 10.5 0
\put {$\bullet$} at 13.5 0
\put {$\bigstar$} at 12 -1

\setlinear
\plot 12 -1 10.5 0 /
\plot 12 -1 13.5 0 /
\plot 10.5 0 9.5 1 /
\plot 10.5 0 10.5 1 /
\plot 10.5 0 11.5 1 /
\plot 13.5 0 12.5 1 /
\plot 13.5 0 13.5 1 /
\plot 13.5 0 14.5 1 /

\endpicture \]

\noindent but also other 3-tuples of trees still of in-degree type $A \, A \, A$ may 
correspond to permutative automorphisms of $\O_3$, e.g.

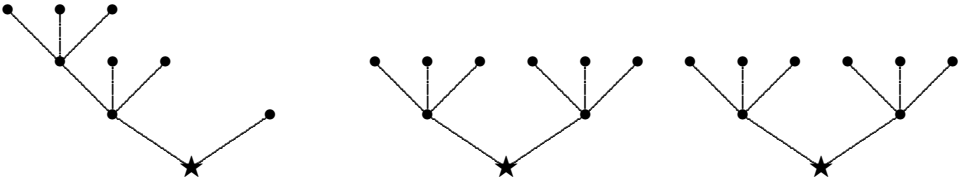
\vspace{-2mm}
\[ \beginpicture
\setcoordinatesystem units <0.7cm,0.7cm>
\setplotarea x from 10 to 15, y from -1 to 1

\put {$\bullet$} at -1.5 2
\put {$\bullet$} at -2.5 2
\put {$\bullet$} at -3.5 2
\put {$\bullet$} at -0.5 1
\put {$\bullet$} at -1.5 1
\put {$\bullet$} at -2.5 1
\put {$\bullet$} at -1.5 0
\put {$\bullet$} at 1.5 0
\put {$\bigstar$} at 0 -1

\setlinear
\plot 0 -1 -1.5 0 /
\plot 0 -1 1.5 0 /
\plot -1.5 0 -2.5 1 /
\plot -1.5 0 -1.5 1 /
\plot -1.5 0 -0.5 1 /
\plot -2.5 1 -3.5 2 /
\plot -2.5 1 -2.5 2 /
\plot -2.5 1 -1.5 2 /

\put {$\bullet$} at 6.5 1
\put {$\bullet$} at 7.5 1
\put {$\bullet$} at 8.5 1
\put {$\bullet$} at 3.5 1
\put {$\bullet$} at 4.5 1
\put {$\bullet$} at 5.5 1
\put {$\bullet$} at 4.5 0
\put {$\bullet$} at 7.5 0
\put {$\bigstar$} at 6 -1

\setlinear
\plot 6 -1 4.5 0 /
\plot 6 -1 7.5 0 /
\plot 4.5 0 3.5 1 /
\plot 4.5 0 4.5 1 /
\plot 4.5 0 5.5 1 /
\plot 7.5 0 6.5 1 /
\plot 7.5 0 7.5 1 /
\plot 7.5 0 8.5 1 /

\put {$\bullet$} at 12.5 1
\put {$\bullet$} at 13.5 1
\put {$\bullet$} at 14.5 1
\put {$\bullet$} at 9.5 1
\put {$\bullet$} at 10.5 1
\put {$\bullet$} at 11.5 1
\put {$\bullet$} at 10.5 0
\put {$\bullet$} at 13.5 0
\put {$\bigstar$} at 12 -1

\setlinear
\plot 12 -1 10.5 0 /
\plot 12 -1 13.5 0 /
\plot 10.5 0 9.5 1 /
\plot 10.5 0 10.5 1 /
\plot 10.5 0 11.5 1 /
\plot 13.5 0 12.5 1 /
\plot 13.5 0 13.5 1 /
\plot 13.5 0 14.5 1 /

\endpicture \]
In fact, type $A$ is comprised of the two kinds of trees entering the above two 3-tuples. 

As deduced in \cite{CKS2} after very long and tedious computer-assisted computations,  
we can report that,
among the permutations already selected on the basis of condition $(b)$,
the total number of permutations for $n=k=3$ satisfying condition
$(d)$ is  
$$907\:044\cdot9!=329\:148\:126\:720 \ . $$
This result relies very much on the extensive set of datas already collected in \cite{CKS1}.
For each of the $7\:390$ representatives in the set $\F$ 
we found the induced    
permutations that satisfy condition $(d)$;
it took about 7 processor years to compute.

In Table \ref{some_sat_condD}, 
we indicate how many instances of each
in-degree type 
have some permutations satisfying condition $(d)$
and how many instances have none.
In Table \ref{reduced_condD}, the 7 390 representative tree tuples are 
counted (second column headed $\#f$)
according to exactly how many induced permutations satisfy condition $(d)$
(first column headed $\#\sigma$),
and according to the combined relabelling and repositioning orbit size   
(third column headed $\#o$).
The fourth entry in each row is the product of the first three entries;
so the sum of the fourth column is the given figure.


\begin{table}
\caption{In-degree types of permutations satisfying condition (d)}
\label{some_sat_condD}
\begin{center}
\begin{tabular}{|rrr|r|r|r|}
\hline
\multicolumn{3}{|c|}{\small ID types} & some & none & total \\
\hline
 A &  A &  A &  290 & 1 878 & 2 168 \\
 A &  B &  B &  611 & 2 171 & 2 782 \\
 A &  C &  D &  86 & 864 & 950 \\
 A &  E &  E &  290 & 782 & 1 072 \\ 
 A &  F &  G &  35 & 357 & 392 \\
 A &  H &  H &  12 & 14 & 26 \\
\hline
\multicolumn{3}{|c|}{total} & 1 324  & 6 066  & 7 390 \\
\hline
\end{tabular}
\end{center}
\end{table}

\begin{table}[h]
\caption{Number of permutations $\sigma$ satisfying condition (d) per $f$.}
\label{reduced_condD}
\begin{center}
\begin{small}
\begin{tabular}{|r@{\ \ \ \ \ \ }r@{\ \ \ \ \ \ \ \ }r@{\ \ \ \ }r|}
\hline
$\#\sigma$ & $\#f$  & $\#o$ & $\#\sigma\cdot\#f\cdot\#o$\\
\hline
0 & 6 066  & $6\cdot9!$ & $0\cdot9!$ \\
24 & 22  & $6\cdot9!$ & $3\:168\cdot9!$ \\  
48 & 288 & $6\cdot9!$ & $82\:944\cdot9!$ \\  
60 & 9 & $6\cdot9!$ & $3\:240\cdot9!$ \\   
72 & 10  & $6\cdot9!$ & $4\:320\cdot9!$ \\
84 & 47 & $6\cdot9!$ & $23\:688\cdot9!$ \\ 
96 & 213  &  $6\cdot9!$ & $122\:688\cdot9!$ \\
96 &   6  &  $3\cdot9!$ & $1\:728\cdot9!$ \\
108 & 103 & $6\cdot9!$ & $66\:744\cdot9!$ \\
120 & 74 & $6\cdot9!$ & $53\:280\cdot9!$ \\
132 & 107 & $6\cdot9!$ & $84\:744\cdot9!$ \\
144 & 111 & $6\cdot9!$ & $95\:904\cdot9!$ \\
156 & 121 & $6\cdot9!$ & $113\:256\cdot9!$ \\
168 & 23 & $6\cdot9!$ & $23\:184\cdot9!$ \\
180 & 3  & $3\cdot9!$ & $1\:620\cdot9!$ \\
192 & 57 & $6\cdot9!$ & $65\:664\cdot9!$ \\
192 &  8 & $3\cdot9!$ & $4\:608\cdot9!$ \\
204 & 26 & $6\cdot9!$ & $31\:824\cdot9!$ \\ 
204 &  4 & $3\cdot9!$ & $2\:448\cdot9!$ \\  
216 & 11 & $6\cdot9!$ & $14\:256\cdot9!$ \\
216 &  7 & $3\cdot9!$ & $4\:536\cdot9!$ \\  
228 & 27 & $6\cdot9!$ & $36\:936\cdot9!$ \\ 
240 & 38 & $6\cdot9!$ & $54\:720\cdot9!$ \\  
312 & 4 &  $6\cdot9!$ & $7\:488\cdot9!$ \\ 
312 & 4 & $3\cdot9!$  & $3\:744\cdot9!$ \\
312 &  1  & $1\cdot9!$& $312\cdot9!$ \\
\hline
    & 7 390 & & $907\:044\cdot9!$ \\
\hline
\end{tabular}
\end{small}
\end{center}
\end{table}

All in all, taking into account the results in \cite{CS,CKS1},
Table \ref{updated} summarizes the up-to-date enumeration of permutations
providing automorphisms of $\O_n$ and 
(in brackets, in the second line) 
of $\D_n$:


\begin{table}
\caption{Number of permutative automorphisms of $\O_n$ (and of $\D_n$) at level $k$}
\label{updated}
\begin{center}
\begin{tabular}{|l||l|l|l|l|}
\hline
$n \setminus k$ & 1 & 2 & 3 & 4\\
\hline\hline
2 & 2 &  4 & 48 & 564,480 \\
  & (2) &  (8) & (324) & (175,472,640) \\
\hline
3 & 6 &  576 & 329,148,126,720 & \\
  & (6) & (5184) & (161,536,753,300,930,560) & \\
\hline
4 & 24 & 5,771,520 & & \\
  & (24) & (1,791,590,400) & & \\
\hline
\end{tabular}
\end{center}
\end{table}
\medskip

\begin{problem}
Extend the results summarized in Tables \ref{some_sat_condD}, \ref{reduced_condD}, and \ref{updated} 
to include a wider range of parameters, possibly developing new computational techniques to this end.
\end{problem}

\section{The restricted Weyl group of $\O_n$}\label{rweyl}

We recall from \cite{Cun2} that $\auto(\O_n,\D_n)$ is the normalizer of $\auto_{\D_n}(\O_n)$ in 
$\auto(\O_n)$ and it can be also described as the group $\lambda(\N_{\O_n}(\D_n))^{-1}$ of 
automorphisms of $\O_n$ induced by elements in the normalizer $\N_{\O_n}(\D_n)$.
Furthermore, using \cite{P}, one can show that $\auto(\O_n,\D_n)$ has the structure of a 
semidirect product $\auto_{\D_n}(\O_n) \rtimes \lambda(\SS_n)^{-1}$ \cite{CS}.
In particular, the group $\lambda(\P_n)^{-1}$ is isomorphic with the quotient of 
the group $\auto(\O_n,\D_n)\cap\auto(\O_n,\F_n)$ by its normal subgroup $\auto_{\D_n}(\O_n)$. 
We call it the {\em restricted Weyl group} of $\O_n$, cf. \cite{Cun2,CS}. 
We also note that every unital endomorphism 
of $\O_n$ which fixes the diagonal $\D_n$ point-wise is automatically surjective, i.e. it is an 
element of $\auto_{\D_n}(\O_n)$ \cite[Proposition 3.2]{Conti} 
and that it is easy to construct product-type automorphisms of $\D_n$ that do not extend to (possibly proper) endomorphisms of $\O_n$ \cite[Proposition 3.1]{Conti}. 
A simple example of such an automorphism of $\D_2$ is given by $\otimes_{i=1}^\infty {\rm Ad}(u_i)$, where $u_i = 1$
for $i$ even and $u_i =( \begin{smallmatrix} 0 & 1 \\ 1 & 0 \end{smallmatrix})$ for $i$ odd 
and we have realized $\D_2$ as an infinite tensor product over $\mathbb N$ of diagonal matrices of size $2$. 
In particular, it becomes important to characterize those automorphisms of $\D_n$ that can be obtained by restricting automorphisms (or even endomorphisms) of $\O_n$. 
As a variation on the theme, we mention the following

\begin{problem} Find necessary and sufficient conditions for an automorphism of $\D_n$ to extend to an automorphism or a proper endomorphism of $\F_n$, respectively.
\end{problem}

In \cite{CHS}, a subgroup $\fG$ of $\auto(\D_n)$ was defined. It consists of those automorphisms 
$\alpha$ for which there exists an $m$ such that both $\alpha\varphi^m$ and $\alpha^{-1}\varphi^m$ 
commute with the shift $\varphi$. 

\begin{theorem}[\cite{CHS}]\label{weyl}
The restriction $r:\lambda(\P_n)^{-1}\to\fG$ is a group isomorphism. 
\end{theorem}

Recall that the spectrum of $\D_n$ may be naturally identified with the full one-sided 
$n$-shift space $X_n$. 
The above theorem identifies the restricted Weyl group of $\O_n$ with the group of those 
homeomorphisms of $X_n$ which together with their inverses eventually commute with the shift. 
In a sense, this provides an answer to a question raised by Cuntz in \cite{Cun2}. 

\vspace{2mm}
We denote $\fIG=\{\Ad(u)|_{\D_n}:u\in\P_n\}$. This is a normal subgroup of $\fG$, 
since for $u\in\P_n^k$ we have $\Ad(u)\varphi^k=\varphi^k$. We also denote 
by $\operatorname{Inn}\lambda(\P_n)^{-1}$ the normal subgroup of $\lambda(\P_n)^{-1}$ 
consisting of all inner permutative automorphisms $\{\Ad(u):u\in\P_n\}$. We call the quotient 
$\lambda(\P_n)^{-1}/\operatorname{Inn}\lambda(\P_n)^{-1}$ the {\em restricted outer Weyl group} 
of $\O_n$. It follows from Theorem \ref{weyl} that the restricted outer Weyl group of 
$\O_n$ is naturally isomorphic to the quotient $\fG/\fIG$. Further analysis reveals that this group in turn is 
related to automorphisms of the two-sided shift. Indeed, 
let $\auto(\Sigma_n)$ denote the group of automorphisms of the full two-sided $n$-shift
(that is, the group of homeomorphisms of the full two-sided $n$-shift space $\Sigma_n$ 
that commute with the two-sided shift $\sigma$)
and let $\la \sigma \ra$ be its subgroup generated by the two-sided shift $\sigma$. 
It is known that $\la \sigma \ra$ coincides with the center of $\auto(\Sigma_n)$. 

\begin{theorem}[\cite{CHS}]\label{outerweyl}
There is a natural embedding of the group $\lambda(\P_n)^{-1}/\operatorname{Inn}
\lambda(\P_n)^{-1}$ into $\auto(\Sigma_n)/\la \sigma \ra$. If $n$ is prime then this embedding 
is surjective and thus the two groups are isomorphic. 
\end{theorem}

The above theorem establishes a useful correspondences between permutative automorphisms of 
the Cuntz algebra $\O_n$ and automorphisms of a classical dynamical system. It opens up 
very attractive possibilities for two-fold applications: of topological dynamics to the study of 
automorphisms of a simple, purely infinite $C^*$-algebra, and of algebraic methods available for 
$\O_n$ to the study of symbolic dynamical systems. Thanks to a combined effort of a number 
of researchers (see \cite{Kit} and \cite{LM}) several interesting properties of the group 
$\auto(\Sigma_n)/\la \sigma \ra$ are known: it is countable, residually finite, contains all finite groups, and 
contains all free products of finitely many cyclic groups. However, a number of questions remain to date 
unsolved (see \cite{Boy}). For example, is it generated by elements of finite order? 
And most importantly, is $\auto(\Sigma_n)/\la \sigma \ra$ isomorphic to $\auto(\Sigma_m)/\la \sigma \ra$ 
(as an abstract group) when $n\neq m$ are prime?

\vspace{2mm}
Along with automorphisms of the two-sided shift, automorphisms of the full one-sided shift $X_n$ 
have been extensively studied (see \cite{Kit} and \cite{LM}). As shown in \cite{CHS}, each 
element of $\auto(X_n)$ (viewed as an element of $\auto(\D_n)$) admits an extension to 
an outer permutative automorphism of $\O_n$. This leads to the following. 

\begin{theorem}\label{extendautos}
There exists a natural (given by extensions) embedding of $\auto(X_n)$ into 
$\lambda(\P_n)^{-1}/\operatorname{Inn}\lambda(\P_n)^{-1}$. 
\end{theorem}

\section{Proper endomorphisms}\label{propend}

In this section, we mainly deal with proper endomorphisms of $\O_n$ which globally 
preserve either the core UHF-subalgebra $\F_n$ or the diagonal MASA $\D_n$. Two main 
references for the results reviewed below are \cite{CRS} and \cite{HSS}. 

\subsection{Endomorphisms preserving $\F_n$}\label{endpresUHF}

Cuntz showed in \cite{Cun2} that if a unitary $w$ belongs to $\F_n$ 
then the corresponding endomorphisms $\lambda_w$ globally preserves 
$\F_n$. The reversed implication was left open in \cite{Cun2}. This question was 
finally answered to the negative in \cite{HSS}, where a number of counterexamples 
were produced. The main method for finding such counterexamples is the following. 

Let $u$ be a unitary in $\O_n$ and let $v$ be a unitary in the relative commutant
$\lambda_u(\F_n)' \cap \O_n$. Then the three endomorphisms
$\lambda_u$, $\lambda_{vu}$, and $\lambda_{u\varphi(v)}$ coincide on  $\F_n$. 
Assume further that $u\in\F_n$, and let $w$ equal either $vu$ or $u\varphi(v)$. Then
$\lambda_u(\F_n)\subseteq\F_n$ and thus $\lambda_w(\F_n)\subseteq\F_n$.
However, $w$ belongs to $\F_n$ if and only if $v$ does.

The above observation shows how to construct examples of unitaries $w$ outside $\F_n$
for which nevertheless $\lambda_w(\F_n) \subseteq \F_n$. To this end, it suffices to find a unitary
$u\in\F_n$ such that the relative commutant $\lambda_u(\F_n)' \cap \O_n$ is not contained in
$\F_n$. This is possible. In fact, one can even find unitaries in a matrix algebra $\F_n^k$
such that $\lambda_u(\O_n)'\cap\O_n$ is not contained in $\F_n$. The existence of such unitaries
was demonstrated in \cite{CP}. The relative commutant $\lambda_u(\O_n)'\cap\O_n$
coincides with the space $(\lambda_u,\lambda_u)$ of self-intertwiners of the
endomorphism $\lambda_u$, which can be computed as
\begin{equation} \label{relcomm}
(\lambda_u,\lambda_u) = \{ x\in\O_n:x=({\rm Ad} u\circ\varphi)(x)\} \ .
\end{equation}
In \cite{CRS}, an explicit example was given of a permutation unitary $u\in\P_2^4$ and a 
unitary $v$ in $\SS_2\setminus\P_2$ such that $v\in(\lambda_u,\lambda_u)$.
Notice that $\lambda_{u\varphi(v)}(\F_n) = \lambda_u(\F_n)$ naturally 
gives rise to a subfactor of the A.F.D. $II_1$ factor with finite Jones index.

\begin{problem}
Provide a combinatorial algorithm to construct and possibly ``classify'' 
pairs $(u,v)$ with $u \in \P_n^k$ and $v \in (\lambda_u,\lambda_u) \cap (\SS_n\setminus \P_n)$.
\end{problem}

An alert reader could spot intriguing resemblance of this problem 
with the classification of the so-called modular invariants (see \cite{BE}), although it 
is possible that this is nothing more than a formal analogy.
%

Furthermore, in \cite{CRS} a striking example was found of a unitary element $u \in \F_2$ for which 
the relative commutant  $\lambda_u(\O_2)' \cap \O_2$ contains a unital copy of $\O_2$. In this 
case the proof is non-constructive and involves a modification of R{\o}rdam's proof of the 
isomorphism $\O_2\otimes\O_2\cong\O_2$, \cite{R}. As a corollary, one obtains 
existence of a unital $*$-homomorphism $\sigma \colon \O_2 \otimes \O_2 \to \O_2$ such that 
$\sigma(\F_2 \otimes   \F_2) \subseteq \F_2$. It is not clear though whether such a $\sigma$ 
can be an isomorphism. 

\begin{problem}
Does there exist an isomorphism $\sigma \colon \O_2 \otimes \O_2 \to \O_2$ such that 
$\sigma(\F_2 \otimes   \F_2) \subseteq \F_2$ or, better yet, $\sigma(\F_2 \otimes   \F_2)= \F_2$?
\end{problem}

At present, we still do not know whether the above described method captures all possible cases or not,  
and thus we would like to pose the following problem. 

\begin{problem}
Does there exist a unitary $w\in\O_n$ such that $\lambda_w(\F_n)\subseteq\F_n$ but 
there is no unitary $u\in\F_n$ such that $\lambda_w|_{\F_n}=\lambda_u|_{\F_n}$?
\end{problem}

Under certain additional assumptions, condition $\lambda_w(\F_n)\subseteq\F_n$ implies $w\in\F_n$, 
\cite{CRS}. In particular, this happens when: 
\begin{description}
\item[\hspace{2mm}{\rm (i)}]  $\lambda_w(\F_n) = \F_n$. If moreover $\lambda_w|_{\F_n}=\id$ then 
$w=t1$, $t\in U(1)$, and thus $\lambda_w$ is a gauge automorphism of $\O_n$;
\item[\hspace{1mm}{\rm (ii)}]  $\lambda_w \in {\rm Aut}(\O_n)$;
\item[{\rm (iii)}] $\lambda_w(\F_n)'\cap\O_n=\bC 1$;
\item[{\rm (iv)}] $w \in \SS_n$ and $\D_n \subseteq \lambda_w(\F_n)$. 
\end{description}

\subsection{Endomorphisms preserving $\D_n$}\label{endpresdiag}

Cuntz showed in \cite{Cun2} that if a unitary $w$ belongs to the normalizer $\N_{\O_n}(\D_n)$ 
of the diagonal $\D_n$ in $\O_n$ then the corresponding endomorphism $\lambda_w$ globally preserves 
$\D_n$. The reversed implication was left open in \cite{Cun2}. This problem was investigated in depth 
in \cite{HSS}. In particular, examples of unitaries $w\not\in\N_{\O_n}(\D_n)$ such that $\lambda_w(\D_n)
\subseteq\D_n$ were found therein, and the following convenient criterion of global preservation 
of $\D_n$ was given. 

\begin{theorem}\label{finitecond}
Let $k \in \bN$ and let $w\in\U(\F_n^k)$. For $i,j=1,\ldots,n$ let $E_{ij}:\F_n^k \rightarrow 
\F_n^{k-1}$ be linear maps determined by the condition that 
$a=\sum_{i,j=1}^n E_{ij}(a)\varphi^{k-1}(S_iS_j^*)$ for all $a\in\F_n^k$. Define 
by induction an increasing sequence of unital selfadjoint subspaces $\mathfrak{W}_r$ of $\F_n^{k-1}$ so that
\[ \begin{aligned}
\mathfrak{S}_1 & = \sp\{ E_{jj}(wxw^*): x\in\D_n^1,\; j=1,\ldots,n \}, \\
\widetilde{\mathfrak{S}}_{r+1} & = \sp\{ E_{jj}((\Ad w\circ\varphi)(x)):
x\in\mathfrak{S}_r,\; j=1,\ldots,n\}, \\
\mathfrak{S}_{r+1} & = \mathfrak{S}_r + \widetilde{\mathfrak{S}}_{r+1}.
\end{aligned} \]
We agree that $\mathfrak{S}_0={\mathbb C}1$.
Let $R$ be the smallest integer such that $\mathfrak{S}_R=\mathfrak{S}_{R-1}$. Then $\lambda_w(\D_n)
\subseteq\D_n$ if and only if $\lambda_w(\D_n^R)\subseteq\D_n$.
\end{theorem}

The above theorem leads to the following corollary, \cite{HSS}. 

\begin{corollary}\label{easyinv}
Let $w$ be a unitary in $\F_n^k$. If $w\D_n^1w^* =\varphi^{k-1}(\D_n^1)$
then $\lambda_w(\D_n)\subseteq\D_n$. Thus if $u\in\F_n^k$, $z\in\U(\F_n^1)$ and 
$u(z\D_n^1z^*)u^* =\varphi^{k-1}(z\D_n^1z^*)$
then $\A=\lambda_z(\D_n)$ is $\lambda_u$-invariant.
\end{corollary}

The second part of the above corollary deals with one of the motivations for investigations of the 
question when $\lambda_w$ preserves $\D_n$. Namely, this information can be useful when 
searching for MASAs of $\O_n$ globally invariant under an endomorphism. The simplest examples 
involve product type standard MASAs, arising as $\lambda_z(\D_n)$ for some 
Bogolubov automorphism $\lambda_z$ of $\O_n$, $z\in\U(\F_n^1)$. Existence of invariant MASAs 
is in turn helpful in determining entropy of an endomorphism, as demonstrated in \cite{SZ,Sk}. 

\bigskip\noindent 
{\bf Acknowledgements.}
The work of J. H. Hong was supported by National Research Foundation of Korea
Grant funded by the Korean Government (KRF--2008--313-C00039).
The work of W. Szyma{\'n}ski was supported by: the FNU Rammebevilling grant `Operator
algebras and applications' (2009--2011), the Marie Curie Research Training Network
MRTN-CT-2006-031962 EU-NCG, the NordForsk Research Network `Operator algebra and dynamics', 
and the EPSRC Grant EP/I002316/1.

%
%
%
%

\end{document}